\documentclass[10pt, a4paper]{article}
\usepackage {fullpage}						
\usepackage {ucs}
\usepackage [utf8x]{inputenc}
\usepackage {graphicx}
\usepackage [bookmarks=false, pdfstartview={FitH}, colorlinks=true, linkcolor=blue, citecolor=blue, urlcolor=blue]{hyperref} 

\usepackage {amsfonts,amsmath,amsthm,amssymb}
\usepackage {latexsym, bbm}
\usepackage {graphics, graphicx}
\usepackage {fancyhdr}
\usepackage {epstopdf} 						

\usepackage {multicol}						
\usepackage {multirow}  					
\usepackage {soul} 							

\usepackage {authblk} 						
\usepackage {datetime} 						
\usepackage {subfigure}						
\usepackage {parskip}						
\usepackage {setspace}						
\usepackage[font=small,textfont=it,width=0.97\textwidth]{caption}		
\usepackage {mathtools}						
\usepackage {array}

\usepackage{subeqnarray}	
\usepackage{overpic}		

\usepackage[usenames,dvipsnames,svgnames,table]{xcolor}

\renewcommand {\(} {\left(}
\renewcommand {\)} {\right)}

\newcommand {\R} {{\rm I\!R}}
\newcommand {\N} {{\rm I\!N}}

\newcommand{\w}{\mathbf w}
\newcommand{\rmT}{{\rm T}}

\def\subtext#1{{\scriptsize\mbox{#1}}}

\theoremstyle{definition}
\newtheorem{experiment}{Experiment}[section]

\begin{document}

\title{Numerical simulation of a contractivity based  multiscale cancer invasion model}
\author{
Niklas Kolbe\thanks{Institute of Mathematics, Johannes Gutenberg-University, Mainz, Germany\hfill {\tt kolbe@uni-mainz.de}}~{},
M\'aria Luk\'a\v{c}ov\'a-Medvid'ov\'a\thanks{Institute of Mathematics, Johannes Gutenberg-University, Mainz, Germany\hfill {\tt lukacova@uni-mainz.de}}~{},
Nikolaos Sfakianakis\thanks{Institute of Mathematics, Johannes Gutenberg-University, Mainz, Germany\hfill {\tt sfakiana@uni-mainz.de}}~{},
Bettina Wiebe\thanks{Institute of Mathematics, Johannes Gutenberg-University, Mainz, Germany\hfill {\tt b.wiebe@uni-mainz.de}}~{}
}

\maketitle

\abstract{We present a problem-suited numerical method for a particularly challenging cancer invasion model. This model is a multiscale haptotaxis advection-reaction-diffusion system that describes the macroscopic dynamics of two types of cancer cells coupled with microscopic dynamics of the cells adhesion on the extracellular matrix. The difficulties to overcome arises from the non-constant advection and diffusion coefficients, a time delay term, as well as stiff reaction terms.
\\
Our numerical method is a second order finite volume implicit-explicit scheme adjusted to include a) non-constant diffusion coefficients in the implicit part, b) an interpolation technique for the time delay, and c) a restriction on the time increment for the stiff reaction terms.
}
\section{Introduction} \sectionmark{Introduction}
The primer objectives in cancer research are to understand the causes of cancer in order to develop strategies for its diagnosis and treatment. The overall effort involves the medical science, biology, chemistry, physics, computer science, and mathematics. The contribution of mathematics, in particular, spans from the modelling of the relevant biological processes, to the analysis of the developed models, and their numerical simualations. The range though of applications of the mathematical models covers a wide range of processes from intracellular bio-chemical reactions to cancer growth, its metastasis and treatment, e.g. \cite{Nordling.1953, Armitage.1954, Fisher.1958,Lukacova.2012,Alt.1985,Preziosi.2003,Perumpanani.1996, Anderson.2000, Gerisch.2008, Szymanska.2009, Painter.2011, Ganguly.2006, Michor.2008, Czochra.2012, 7, Johnston.2010, Vainstein.2012, Gupta.2009}. 
 
In this work we focus in the first step of cancer metastasis ---and one of the ``hallmarks of cancer''--- the invasion of the \textit{extracellular matrix}\index{Extracellular Matrix} (ECM). Our study involves the existence of a secondary group of cancer cells within the main body of the tumour that exhibits \textit{stem-cell-like}\index{cell!stem} properties. This secondary group of cancer cells seem to stem from the ``original'' cancer cells via a cellular differentiation program that can be found also in normal tissue, the \textit{Epithelial-Mesenchymal Transition}\index{Epithelial-Mesenchymal Transition} (EMT). Both the EMT and its reverse process, \textit{Mesenchymal-Epithelial Transition} (MET) participate in several developmental processes including embryogenesis, wound healing, and fibrosis, \cite{Thiery.2002, Mani.2008, Katsuno.2013, Singh.2010, 2}.

The two types of cancer cells possess different cell \textit{proliferation}\index{cell!proliferation} rates and motility properties, and present different levels of cellular \textit{potency}\index{cell!potency}. The secondary group, in particular, exhibits lower (if any) proliferation rates, stem cell-like properties such as self-renewal and cellular differentiation. These cells are more resilient to cancer therapies and they are able to metastasize. While the bulk consists mostly of the ``original'' cancer cells, the secondary family constitutes the smaller part of the tumour, \cite{Gupta.2009, Reya.2001}.
 
The motility mechanism of the cancer cells responds to alterations and gradients in the chemical environment of the tumour (a process termed \textit{chemotaxis}\index{chemotaxis}), and in the cellular adhesions sites located on the ECM (a process termed \textit{haptotaxis}\index{haptotaxis}). From a mathematical perspective the study of several forms of -taxis has been an active research field in the last decades. The derived models are typically \textit{Keller-Segel}\index{Keller-Segel} (KS) type systems \cite{KS.1970, Patlak.1953}, where the participating quantities are described macroscopically in the sense of densities. By including also interactions between the cancer cells and the extracellular environment the resulting models take the form of \textit{Advection-Reaction-Diffusion}\index{Advection-Reaction-Diffusion} (ARD) systems, see e.g. \cite{Alt.1985, Bellomo.2008, Anderson.2000, Andasari.2011, Domschke.2014, Stinner.2015, Sfakianakis.2015, Sfakianakis.2016}.

The solutions of these models exhibit typically complex dynamical behavior manifested in the form of merging/emerging concentration or in the form of complex wave phenomena \cite{Chaplain.2005, Domschke.2014, Sfakianakis.2016}. Moreover, since these models are close to the classical KS systems, the possibility of a blow-up ---if not analytically excluded--- should be numerically investigated. Due to such dynamical behaviour these dynamics special and problem specific numerical treatments are needed, \cite{Gerisch.2008, Sfakianakis.2014b, Sfakianakis.2015}.

In the current paper our aim is to contribute in this direction by presenting our problem-suited method for a particular ARD cancer invasion haptotaxis model that was proposed in \cite{Stinner.2015}. This model features several numerically challenging properties: non constant advection and diffusion coefficients, non-local time delay, and stiff reaction terms, see (\ref{F:eq:orig2}).  

The rest of the paper is structured as follows: in Section \ref{F:sec:model} we present and discuss briefly the cancer invasion model. In Section \ref{F:sec:method} we address the numerical method we employ, comment on its properties, and on the special treatment of its terms. In Section \ref{F:sec:results} we present our numerical findings and discuss their implication in terms of the model. 

\section{Mathematical model}\label{F:sec:model} \sectionmark{Model}
The model we investigate is a cancer invasion ARD system of the KS spirit that primarily features two families of cancer cells and includes \textit{contractivity}: a measure of the strength of the cell-matrix adhesions, see \cite{Stinner.2015} In some more details, the following properties are assumed by the model:
\begin{itemize}
	\item[--] The ``original'' cancer cells (henceforth \textit{proliferative}) proliferate and do not migrate or otherwise translocate. The stem-like cancer cells (henceforth \textit{migratory}) migrate but do not proliferate. 
	\item[--] The motility mechanism of the migratory cancer cells responds to the (possibly) non-uniform distribution of adhesion sites located on the ECM. The induced haptotactic movement is modelled by a combination of advection and  diffusion. 
	\item[--] There exist a \textit{bidirectional} transition between the two families of cancer cells, modelling the parallel action of EMT and MET. Both are assumed to take place with constant rates. 
	\item[--] The ECM is a dynamic structure that is degraded by the cancer cells, and constantly remodelled. Remodelling is self-induced, i.e. new ECM is sprout from existing ECM. 
	\item[--] The proliferation of the cancer cells and the remodelling of the ECM is limited by the locally available space. This effect is modelled by a {volume filling} term. 
	\item[--] Both the diffusion and the advection of the migratory cancer cells are governed by non-uniform coefficients depending on the contractivity. This property determines the ``strength'' of the cell \textit{migration}\index{cell!migration}. 
	\item[--] The contractivity depends ---in a time delay way--- on the local amount of the ECM-bound \textit{integrins}\index{integrins}, which are attached on the ECM. They degrade with a constant rate and are reproduced with a preferable maximum local density. 
	\item[--] The dynamics of the integrins and the {contractivity} take place on a \textit{microscopic time scale} \index{time scale!microscopic} that is faster than the \textit{macroscopic time scale}\index{time scale!macroscopic} of the dynamics of the cancer cells.
\end{itemize}
Altogether the model reads:
\renewcommand{\arraystretch}{1.5}
\begin{equation}\label{F:eq:orig}
\left\{\begin{array}{*1{>{\displaystyle}r}c*1{>{\displaystyle}l}}
		\partial_t c_1&=&\mu_c c_1(1-(c_1+c_2)-\eta_1v)+\gamma c_2-\lambda c_1 \\
		\partial_t c_2&=&\nabla\cdot \(D_c \frac{\kappa}{1+(c_1+c_2)v}\nabla c_2 \)-\nabla\cdot\(
			D_h\frac{ \kappa v}{1+v}c_2\nabla v \)+\lambda c_1-\gamma c_2\\
		\partial_t v&=&-\delta_v(c_1+c_2)v+\mu_v v(1- \eta_2(c_1+c_2)-v)\\
		\partial_\vartheta y&=&k_1(1-y)v-k_{-1}y\\
		\partial_\vartheta\kappa&=&-q\kappa+My(\vartheta-\tau)
\end{array}\right.\,,
\end{equation}
where $c_1$, $c_2$ denote the densities of the proliferating and the migrating cancer cells respectively. The densities of the ECM and of the ECM-bound integrins are denoted by $v$ and $y$. The {contractivity} is denoted by $\kappa$, and the {microscopic} and {macroscopic} time scales by $\vartheta$, and $t$ respectively. We assume that the time scales are related in the following way:
\begin{equation} \label{F:eq:micro_time}
	\vartheta = \frac{t}{\chi}\,,
\end{equation}
where the rescaling factor $0 < \chi <1$ is a fixed constant. The \textit{time delay}\index{time delay} $\tau>0$ is also assumed to be a constant.

The system (\ref{F:eq:orig}) is endowed with initial and boundary conditions, see Section~\ref{F:sec:exp.desc}.

\paragraph{Rescaled system.}
Using the time scale relation (\ref{F:eq:micro_time}) we can rescale the model (\ref{F:eq:orig}) and obtain the following system using only the {macroscopic} time variable $t$:
\renewcommand{\arraystretch}{1.5}
\begin{equation}\label{F:eq:orig2}
\left\{\begin{array}{*1{>{\displaystyle}r}c*1{>{\displaystyle}l}} 
		\partial_t c_1&=&\mu_c c_1(1-(c_1+c_2)-\eta_1v)+\gamma c_2-\lambda c_1\\
		\partial_t c_2&=&\nabla\cdot \(D_c \frac{\kappa}{1+(c_1+c_2)v}\nabla 	c_2\)-\nabla\cdot\(
			D_h\frac{ \kappa v}{1+v}c_2\nabla v \)+\lambda c_1-\gamma c_2 \\
		\partial_t v&=&-\delta_v(c_1+c_2)v+\mu_v v(1- \eta_2(c_1+c_2)-v) \\
		\partial_t y&=&\frac{k_1}{\chi}(1-y)v-\frac{k_{-1}}{\chi}y \\
		\partial_t\kappa&=&\displaystyle -\frac{q}{\chi}\kappa+\frac{M}{\chi}y(t-\chi\tau)
\end{array}\right.\,.
\end{equation}

\paragraph{Operator Form.}
The system (\ref{F:eq:orig2}) can be written for convenience in a compact operator form as follows:
\begin{equation}
\label{F:opForm2}
	\w_t=D(\w) - A(\w) + R(\w)\,,
\end{equation}
where $\w:\Omega\times \R_+ \rightarrow \R^5$, with $\w=(c_1, c_2,v, y,\kappa)^T$, and $D$, $A$, $R$ represent the diffusion, advection, and reaction operators, respectively:
\renewcommand{\arraystretch}{1.5}
\[\begin{array}{rlc}
	D(\w)&=\(0,\ \nabla\cdot \(D_c \frac{\kappa}{1+(c_1+c_2)v}\nabla c_2\),\ 0,\ 0,\ 0 \)^T\,,\\[0.7em]
	A(\w)&=\(0,\ \nabla\cdot \(D_H \frac{\kappa v}{1+v}c_2\nabla v \)  ,\ 0, \ 0, \ 0\)^T\,,\\[0.7em]
	R(\w)&=\(
		\begin{array}{c}
			\mu_c c_1(1-(c_1+c_2)-\eta_1v)+\gamma c_2-\lambda c_1\\
			\lambda c_1-\gamma c_2\\
			-\delta_v(c_1+c_2)v+\mu_v v(1- \eta_2(c_1+c_2)-v)\\
			\frac{k_1}{\chi}(1-y)v-\frac{k_{-1}}{\chi}y\\
			-\frac{q}{\chi}\kappa+\frac{M}{\chi}y(\vartheta-\tau)
		\end{array}\)\,.
\end{array}\]
Additionally we set 
\begin{equation}\label{F:RimplII}
R_\subtext{impl}(\vec w)=\(0 ,\ 0 ,\ 0 ,\ \frac{k_1}{\chi}(1-y)v-\frac{k_{-1}}{\chi}y ,\ 0 \)^T\,,
\end{equation}
and 
\begin{equation}\label{F:RexplII}
	R_\subtext{expl}(\vec w)=R(\vec w) - R_\subtext{impl}(\vec w)\,.
\end{equation}

\paragraph{Parameters.}
\label{F:paramII}
For the main experiments we consider the following set of parameters that has been adopted by \cite{Stinner.2015}: 
\renewcommand{\arraystretch}{1.1}
\begin{equation}\label{F:params}
	\left\{
	\begin{array}{*4{>{\displaystyle}l}}
	\mu_c=1,		&\qquad \eta_1=0.05, 	&\qquad \gamma=0.055,		&\qquad \lambda=0.152, \\
	D_c=0.01,	&\qquad D_h=10, \\
	\delta_v=5, 	&\qquad \mu_v=0.3, 		&\qquad \eta_2=0.9, \\
	k_1=2, 			&\qquad k_{-1}=0.06,\\
	q=3, 			&\qquad M=2,\\
	\chi= 0.01, 	&\qquad \tau=20.
	\end{array}
	\right.
\end{equation}
These parameters are adjusted in each particular experiment under investigation, see also Section \ref{F:sec:exp.desc}.

\section{Numerical method}\label{F:sec:method} \sectionmark{Numerical method}
We consider a two-dimensional computational domain $\Omega=(a,b)\times(a,b)\subset \R^2$, which will be subdivided into a finite number of 
regular computational cells of size:
\[
	h=(h_1,\;h_2)^T\ \mbox{ where }\  h_1=\frac{b-a}{L},\ h_2=\frac{b-a}{M}\,.
\]
Here $L,M\in\N$ denotes the resolution of the grid along the $x_1$- and $x_2$-directions, respectively. the total number of grid cells is $N=LM$. The cell centers are located at
\begin{eqnarray}
	\vec x_{1,1} &=& \(a+\frac{h_1}{2}\)\vec{e}_1 + \(a+\frac{h_2}{2}\)\vec{e}_2\,,\\
	\vec x_{i,j} &=& \vec x_{1,1} + (i-1)\, h_1 \, \vec{e}_1 + (j-1)\, h_2\, \vec{e}_2\,,
\end{eqnarray}
for $i=1,\dots,L$, $j=1,\dots,M,$ where $\vec e_1, \vec e_2$ are the unit vectors along the $x_1$- and $x_2$-directions, respectively. Consequently, the computational cells are given by
\[
	C_{i,j}=\left\{ \vec x_{i,j} + \left(\lambda_1\,h_1 , \lambda_2\, h_2\right) ,~\lambda_1,\lambda_2 \in \left[-\frac12, \frac12\right) \right\}\,,
	\ i=1,\dots,L,\ j=1,\dots,M\,.
\]
We introduce a single-index notation for the two-dimensional computational cells using the \textit{lexicographical order}, i.e.
\begin{subeqnarray}
	C_{i,j} &\longrightarrow& C_{i+(j-1)L}\,,\\
	\vec x_{i,j} &\longrightarrow& \vec x_{i+(j-1)L}\,,
\end{subeqnarray}
for $i=1,\dots,L$, $j=1,\dots,M$, and inversely
\begin{subeqnarray}
	C_k &\longrightarrow& C_{k-\lfloor \frac{k-1}{L}\rfloor L ,\lfloor\frac{k-1}{L}\rfloor +1}\,,\\
	\vec x_{k} &\longrightarrow& \vec x_{k-\lfloor\frac{k-1}{L}\rfloor L ,\lfloor\frac{k-1}{L}\rfloor +1}\,, 
\end{subeqnarray}
for $k=1,\dots,N$, where $\lfloor ~\rfloor$ is the \textit{Gauss floor function}\index{Gauss floor function}. We denote moreover by $C_{k\pm \vec e_j}$ the neighbouring cell of $C_k$ along the positive (negative) $\vec e_j$ direction ($j=1,2$). Hence for $k=1,\dots ,N$ we have
\begin{subeqnarray}
	~\hspace{-2.5em}C_{k \pm \vec e_1}=  C_{k-\lfloor\frac{k-1}{L}\rfloor  L \pm  1, \lfloor\frac{k-1}{L}\rfloor +1},&& \mbox{for }k\neq 0,1\, \mbox{mod}\ L,\,\ \mbox{respectively}\,,\\
	~\hspace{-2.5em}C_{k \pm \vec e_2}=  C_{k- \lfloor\frac{k-1}{L}\rfloor L , \lfloor\frac{k-1}{L}\rfloor +1 \pm 1},&& \mbox{for }k\leq  L(M-1),\ k\geq  L+1 ,\,\ \mbox{respectively}\,.
\end{subeqnarray}

\subsection{Space Discretization}

The system (\ref{F:eq:orig2}) is discretized in space by a finite volume method. The approximate solution is represented on every computational cell $C_i$ by a piecewise constant function
\begin{equation}\label{F:approx}
	\w_i(t)\approx \frac{1}{|C_i|}\int_{C_i} \w(x,t)~dx\,.
\end{equation}
Moreover, for $\w_h(\cdot)=\left\{\w_i(\cdot)\right\}_{i=1}^N$, we consider the following approximations of the advection, diffusion, and reaction operators:
\renewcommand{\arraystretch}{1.5}
\begin{equation}\label{F:eq:dscrOper}
\left\{
\begin{array}{*1{>{\displaystyle}r}c*1{>{\displaystyle}l}} 
	{\rm A}_i(\w_h(t)) &\approx& \frac{1}{|C_i|}\int_{C_i} A(\w(x,t))\;dx\,, \\
	{\rm D}_i(\w_h(t)) &\approx& \frac{1}{|C_i|}\int_{C_i} D(\w(x,t))\;dx\,,\\
	{\rm R}_i(\w_h(t)) &\approx& \frac{1}{|C_i|}\int_{C_i} R(\w(x,t))\;dx\,.
\end{array}
\right.
\end{equation}
\paragraph{Reaction.}
	The reaction terms are discretized by a direct evaluation of the reaction operator at the cell centers
	\begin{equation}\label{F:eq:dscr_Reac}
		{\rm R}_i(\w_h(t)) = R(\w_i(t))\,.
	\end{equation}
	
\paragraph{Diffusion.} 
	We denote the discrete diffusion coefficient, see also (\ref{F:eq:orig2}), by
	\[\rmT_i(\w_h(t)) = \frac{D_c\, \kappa_i}{1 + (c_{1,i}+ c_{2,i})v_i},\]
	and define the second component of the discrete diffusion operator using \emph{central differences}
	\begin{equation}\label{F:eq:dscr_Adv}
	\begin{array}{*1{>{\displaystyle}r}*1{>{\displaystyle}l}} 
		[{\rm D}_i(\w_h(t))]_2 = \sum_{j=1}^2& \frac{\rmT_{i-\vec e_j}(\w_h(t)) + \rmT_i(\w_h(t)))}{2 \, h_j^2} c_{2, i - \vec e_j} \\
		 &- \frac{ \rmT_{i-\vec e_j}(\w_h(t)) + 2 \,\rmT_i(\w_h(t)) +   \rmT_{i+\vec e_j}(\w_h(t))}{2 \, h_j^2} c_{2, i} \\
		 &+\frac{\rmT_i(\w_h(t)) + \rmT_{i+\vec e_j}(\w_h(t))}{2 \, h_j^2} c_{2, i + \vec e_j} \,, 
	\end{array}
	\end{equation}
	with all remaining components $[{\rm D}_i(\w_h(t))]_j$, $j=1,3,4,5$ being equal to zero.
	
\paragraph{Advection.} 
	The advection term is discretized using the \emph{central upwind flux}, see \cite{ck,Lukacova.2012}, which in the particular case of the system (\ref{F:eq:orig2}) reads as
	\begin{equation}\label{F:eq:dscr_Adv}
		{\rm A}_i(\w_h(t)) = \sum_{j=1}^2 \frac{1}{h_j}\(0,\ {\rm H}_{i+\vec e_j/2}\(\w_h(t)\)-{\rm H}_{i-\vec e_j/2}\(\w_h(t)\),\ 0,\  0,\  0\)^T\,.
	\end{equation}
	The numerical flux ${\rm H}_{i+\vec e_j/2}$ approximates the flux between the computational cells $C_i$ and $C_{i+\vec e_j}$, $j=1,2$:
	\begin{equation} \label{F:numFlux}
		{\rm H}_{i+\vec e_j/2}(\w_h)=
		\left\{\begin{array}{ll}
			{\rm P}_{i+\vec e_j/2}(\w_h)\, c_{2,i+\vec e_j/2}^{+},	&\textnormal{ if }{\rm P}_{i+\vec e_j/2}(\vec w_h)\geq 0\,, \\
			{\rm P}_{i+\vec e_j/2}(\w_h)\, c_{2,i+\vec e_j/2}^{-},	&\textnormal{ if }{\rm P}_{i+\vec e_j/2}(\vec w_h)< 0\,,
		\end{array}\right.
	\end{equation}
	and ${\rm P}_{i+\vec e_j/2}$ represents the local characteristic speeds as:
	\[
		{\rm P}_{i+\vec e_j/2}(\w_h)= \frac{D_h}{2} \left( \frac{\kappa_i \, v_i }{1 + v_i} + \frac{\kappa_{i+\vec e_j}\,  v_{i+\vec e_j}}{1 + v_{i+\vec e_j}} 
		\right)\frac{v_{i+\vec e_j}-v_{i}}{h_j},		
	\]
		for both space directions $j=1,2$. The interface values $c_{2,i+\vec e_j/2}^{\pm}$ are computed by the linear reconstructions
	\begin{subeqnarray}
		c_{2,i+\vec e_j/2}^{-} &=& c_{2,i} + s_{i}^{(j)}\,,\\
		c_{2,i+\vec e_j/2}^{+} &=&c_{2,i+\vec e_j} - s_{i+\vec e_j}^{(j)}\,,
	\end{subeqnarray}
	where the slopes $s_{i}^{(j)}$ are provided by the \textit{monotonized central}\index{monotonized central limiter} (MC) limiter \cite{van1977towards}
	\begin{equation} \label{F:eq:nc}
		s_i^{(j)}={\rm minmod}\(c_{2,i}-c_{2,i-\vec e_j}, \frac{1}{4} (c_{2,i+\vec e_j} -c_{2,i-\vec e_j}), c_{2,i+\vec e_j}-c_{2,i}\).
	\end{equation}
	The \textit{minmod operator} is given by
	\begin{equation}\label{F:eq:minmod}
		{\rm minmod}(v_1,v_2,v_3)=\left\{
		\begin{array}{ll}
			{\rm max}\{v_1,v_2,v_3\},& \mbox{if}\ v_k<0,\ k=1,2,3,\\
			{\rm min}\{v_1,v_2,v_3\},& \mbox{if}\ v_k>0,\ k=1,2,3,\\
			0,& \mbox{otherwise}.\\
		\end{array}\right.
	\end{equation}

Applying the above and (\ref{F:eq:dscrOper}), (\ref{F:eq:dscr_Reac}), (\ref{F:eq:dscr_Adv}), we obtain the system of the Ordinary Differential Equations (ODEs)
\begin{equation}\label{F:eq:semi-dscr}
	\partial_t \w_h -{\rm A}(\w_h)= {\rm R}(\w_h) + {\rm D}(\w_h)\;.
\end{equation}

\subsection{Time Discretization}\label{F:methods}
Let us consider  $\w_h^n$ a numerical approximation of the solution $\w_h(t)$ of (\ref{F:eq:semi-dscr}) at discrete time instances $t_n$, where $t_n=t_{n-1}+\Delta t_n$. For the choice of the time steps $\Delta t_n$ we refer to Section \ref{F:sec:dt}.

\paragraph{IMEX.}
For the time discretization of (\ref{F:eq:semi-dscr}) we employ an \textit{Implicit-Explicit Runge-Kutta}\index{Implicit-Explicit Runge-Kutta} (IMEX) method of 3rd order of accuracy first proposed in \cite{christopher2001additive}. 

A diagonally implicit Runge-Kutta (RK) scheme is applied to the implicit part and an explicit Runge-Kutta scheme to the explicit part. The scheme can be written in the following way:
\begin{equation}\label{F:eq:IMEX3}
\left\{\begin{array}{*1{>{\displaystyle}r}c*1{>{\displaystyle}l}} 
	\vec w_h^{n+1}=\vec w_h^n +\Delta t_n\Big( \sum\limits_{j=1}^{i-1}b_{j}^E ({\rm -A + R_\subtext{expl}})(t_n + c_j^E \Delta t, \vec W_i)
		+\sum\limits_{j=1}^{i}b_{j}^I({\rm D + R_\subtext{impl}})(\vec W_i) \Big),\\
	\vec W_i= \vec w_h^n +\Delta t_n \Big( \sum\limits_{j=1}^{i-1}a_{ij}^E ({\rm -A + R_\subtext{expl}})(t_n + c_j^E \Delta t, \vec W_i)
		+\sum\limits_{j=1}^{i}a_{ij}^I(\rm D + R_\subtext{impl})(\vec W_i) \Big) . 
\end{array}\right.
\end{equation}		
where $ \vec b^E, \vec c^E \in \R^s$, $A^E\in \R^{s\times s}$, $\vec b^I,\vec c^I\in \R^s$ and $A^I\in \R^{s\times s}$ stand for the explicit, and implicit scheme coefficients, respectively. Note, we approximate the advection operator explicitly in time. Using the splitting of the reaction operator according to (\ref{F:RexplII}) and (\ref{F:RimplII}), the reaction terms are computed in both explicitly and implicitly. We finally compute the stages $\vec W_i$ by solving the linear system in the second equation of (\ref{F:eq:IMEX3}) using the iterative biconjugate gradient stabilized \textit{Krylov subspace method}\index{Krylov subspace method} \cite{Krylov.1931, vdVorst.1992}.


The particular four stage ($s=4$) IMEX method that we employ uses the Butcher Tableau \ref{F:tbl:IMEX} and fulfills several stability conditions like A- and L-stability \cite{christopher2001additive}. 

\begin{table}[t]
\begin{center}
\begin{tabular}{c|cccc}
$0$&&&&\\
$\frac{1767732205903}{2027836641118}$&$\frac{1767732205903}{2027836641118}$&&&\\
$\frac{3}{5}$&$\frac{5535828885825}{10492691773637}$&$\frac{788022342437}{10882634858940}$&&\\
$1$&$\frac{6485989280629}{16251701735622}$&$-\frac{4246266847089}{9704473918619}$&$\frac{10755448449292}{10357097424841}$&\\
\hline
&$\frac{1471266399579}{7840856788654}$ & $-\frac{4482444167858}{7529755066697}$ & $\frac{11266239266428}{11593286722821}$ & $\frac{1767732205903}{4055673282236}$ 
\end{tabular}
\end{center}

\begin{center}
\begin{tabular}{c|cccc}
$0$&0&&&\\
$\frac{1767732205903}{2027836641118}$&$\frac{1767732205903}{4055673282236}$&$\frac{1767732205903}{4055673282236}$&&\\
$\frac{3}{5}$&$\frac{2746238789719}{10658868560708}$&$-\frac{640167445237}{6845629431997}$&$\frac{1767732205903}{4055673282236}$&\\
$1$&$\frac{1471266399579}{7840856788654}$&$-\frac{4482444167858}{7529755066697}$&$\frac{11266239266428}{11593286722821}$&$\frac{1767732205903}{4055673282236}$\\
\hline
&$\frac{1471266399579}{7840856788654}$ & $-\frac{4482444167858}{7529755066697}$ & $\frac{11266239266428}{11593286722821}$ & $\frac{1767732205903}{4055673282236}$ 
\end{tabular}
\end{center}
\caption{Butcher tableaux for the explicit (upper) and the implicit (lower) parts of the third order IMEX scheme (\ref{F:eq:IMEX3}), see also \cite{christopher2001additive}}
\label{F:tbl:IMEX}
\end{table}
\subsection{Treatment of the delay term}
Of particular importance for the system (\ref{F:eq:orig2}) is the \textit{time delay}\index{time delay} term $y(t-\chi\tau)$ that appears in the contractivity equation $\kappa$. We have included this term in the explicit part ${\rm R}_\subtext{expl}$ of the implicit-explicit description (\ref{F:RexplII}) of ${\rm R}$. Consequently an approximation of the delayed component is needed in the explicit part of the IMEX method. 

At stage $j$ of the method (\ref{F:eq:IMEX3}) we evaluate the operator ${\rm R}_\subtext{expl}$ at the time instance $\hat t = t_n + c_j^E \Delta t_n$, thus we need to approximate $y(\hat t-\chi\tau)$. We identify the position of $\hat t-\chi\tau$ (recall that $\chi,\tau\geq 0$) and interpolate between the known values of $y_h$. In some more detail: we consider the time instances $t_1^d \leq t_2^d \leq t_n \leq \hat t$, where $t_1^d \leq \hat t-\chi\tau$, and corresponding densities of the integrins $y_h(t_1^d),~ y_h(t_2^d)~, y_h(t_n)~, y_h(\hat t)$. Then the interpolated value of $y$ is given as:
\[
y_h(\hat t-\chi\tau) = \left\{
\begin{array}{l l}
\displaystyle y_h(t_1^d) + \frac{\hat t -\chi\tau - t_{1}^d }{t_{2}^d - t_{1}^d}\( y_h(t_2^d) -y_h(t_1^d)\),  &\quad \textnormal{if } t_1^d \leq \hat t-\chi\tau < t_2^d,\nonumber\\
\displaystyle y_h(t_2^d) + \frac{\hat t -\chi\tau - t_{2}^d }{t_n - t_{2}^d}\( y_h(t_n) -y_h(t_2^d)\), &\quad \textnormal{if } t_2^d \leq \hat t-\chi\tau < t_n,\nonumber\\
\displaystyle y_h(t_n) + \frac{\hat t -\chi\tau - t_n }{\hat t - t_n}\( y_h(\hat t) -y_h(t_n)\), &\quad \textnormal{if } t_n \leq \hat t-\chi\tau < \hat t,\nonumber
\end{array}\right.
\]
where we use the current numerical solution $y_h(t_n) = [\vec w_h^n]_4$, and previously computed $y_h(t_1^d), y_h(t_2^d)$. Further, we make the assumption that $W_j \approx w_h(\hat t)$ and hence employ $y_h(\hat t) = [W_j]_4$ to approximate the integrin density at time instance $\hat t$.

Having computed the time update $w_h^{n+1}$, we check if the delay time at the next time integration step will overshoot $t_2^d$. Thus, if $t_{n+1}- \chi \tau \geq t^d_2$ we update the time instances and the corresponding numerical solutions used for the interpolation by
$$ t_1^d \leftarrow t_2^d, \quad t_2^d \leftarrow t_n, \quad y_h(t_1^d) \leftarrow y_h(t_2^d), \quad y_h(t_2^d) \leftarrow y_h(t_n).$$ 
\begin{figure}[t]
\centering
\begin{overpic}[width=1\textwidth]{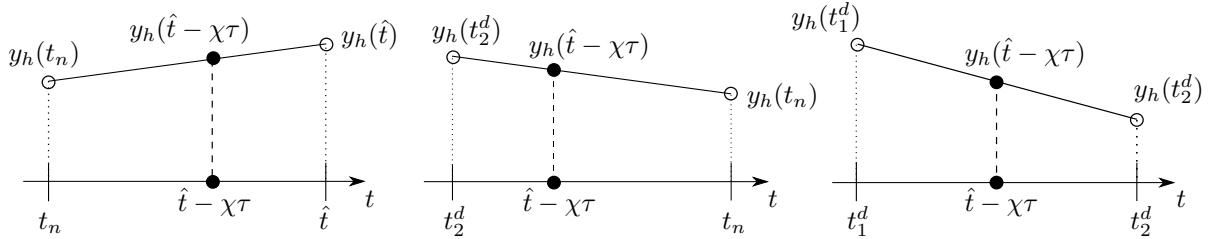}
 \put (30.5,3.3) {$t$}
 \put (63.8,3.3) {$t$}
 \put (97.5,3.3) {$t$}
 \put (15,3) {$\hat t-\chi\tau$}
 \put (11,17.5) {$y_h(\hat t-\chi\tau)$}
 \put (26.7,1.5) {$\hat t$}
 \put (3.5,1.5) {$t_{n}$}
 \put (60,1.5) {$t_{n}$}
 \put (37,1.5) {$t_{2}^d$}
 \put (93.7,1.5) {$t_{2}^d$}
 \put (70.3,1.5) {$t_{1}^d$}
 \put (28.5,17) {$y_h(\hat t)$}
 \put (1,15.5) {$y_h(t_n)$}
 \put (43.5,3) {$\hat t-\chi\tau$}
 \put (44,16.2) {$y_h(\hat t-\chi\tau)$}
 \put (80,3) {$\hat t-\chi\tau$}
 \put (80,15.5) {$y_h(\hat t-\chi\tau)$}
 \put (62,12) {$y_h(t_n)$}
 \put (36,17.4) {$y_h(t^d_2)$}
 \put (65.7,18.5) {$y_h(t^d_1)$}
 \put (94,12.5) {$y_h(t^d_2)$}
\end{overpic}
\caption{Illustration of the three different cases used to interpolate the delay term}
\label{F:fig:interpolate}
\end{figure}

\subsection{Choice of the time step}\label{F:sec:dt}

For stability reasons the time steps are restricted by the characteristic velocities using the CFL condition, \cite{CFL.1928}:
\begin{equation}\label{F:eq:CFL}
	\max_{i=1,\dots, N,\ j=1,2}\; \Delta t_n\, \frac{{\rm P}_{i+\vec e_j/2}}{h_j}\leq 0.5\;.
\end{equation}

Moreover, we note that the ODE subsystem of the last two equations of (\ref{F:eq:orig2}) is stiff due to the large parameter $\frac{1}{\chi}$, cf. (\ref{F:params}). This results in instabilities and inaccuracies in our partly explicit method if the time steps are not further regulated. To cure this problem, we choose $\Delta t_n$ such that the relative change of $\kappa$ remains bounded, i.e.
\begin{equation}\label{F:eq:kappa_dt_restriction}
\frac{\| \kappa_h(t_n) - \tilde \kappa_h(t_n + \Delta t_n)\|_\infty}{\|\kappa_h(t_n) \|_\infty}\leq 0.01\;,
\end{equation}
where $\kappa_h(t_n) = [\vec w_h^n]_5$. Since the time increment $\Delta t_n$ is needed in our method to compute the actual approximate contractivity $\kappa(t_n + \Delta t_n) = [\vec w_h^{n+1}]_5$ we apply (\ref{F:eq:kappa_dt_restriction}) using an estimator $\tilde \kappa_h(t_n + \Delta t_n) \approx \kappa_h(t_n + \Delta t_n)$. It can be seen in (\ref{F:eq:orig2}) that the component $\kappa$ evolves quickly (in physical time) to a quasi-steady-state, in order to keep the second order accuracy, cf also \cite{Wiebe.2016}. After this state is reached, the changes in $\kappa$ are very slow. Hence the restriction (\ref{F:eq:kappa_dt_restriction}) affects the employed time increment only at the beginning of the computation.

The choice of the threshold value $0.01$ in (\ref{F:eq:kappa_dt_restriction}) has followed from numerical experimentation and in order to keep the second order accuracy, see  Tables \ref{F:tbl:L1EOC}--\ref{F:tbl:L2EOC} and Fig. \ref{F:fig:cnv.cpu}.

In practice we compute $\Delta t$ from (\ref{F:eq:CFL}) and (\ref{F:eq:kappa_dt_restriction}) as follows: As a very first step in the IMEX method we compute ${\rm A}(t_n, \vec w_h^n),~ {\rm R}_\subtext{expl}(t_n, \vec w_h^n)$ which are needed for the first stage of the RK updates. In the flux computation we get
\[
a =  \max_{i=1,\dots, N,\ j=1,2}\; \frac{{\rm P}_{i+\vec e_j/2}}{h_j}.
\]
Since the equation for the contractivity includes only reaction terms that are evaluated in the operator $\rm R_\subtext{expl}$, we can employ the forward Euler estimator
$$ \tilde \kappa_h(t_n + \Delta t_n) = \kappa_h(t_n) + \Delta t_n [{\rm R}_\subtext{expl}(t_n, \vec W_i)]_5.$$
We do not actually compute the Euler step $ \tilde \kappa_h(t_n + \Delta t_n)$, merely substitute in (\ref{F:eq:kappa_dt_restriction}) and deduce the time increment
\begin{equation}\label{F:eq:dt_compute}
	\Delta t_n = \min \left\{ \frac{1}{2\, a},~ \frac{\|\kappa_h(t_n) \|_\infty}{100 \, \|[{\rm R}_\subtext{expl}(t_n, \vec W_i)]_5 \|_\infty} \right\},
\end{equation}
before we compute a new time update. In effect, we compute $\Delta t_n$ without placing additional computational burden on the method.

\section{Experimental results} \label{F:sec:results}\sectionmark{Experimental results}

\begin{table}[t] 
  \begin{center}\begin{tabular}{|c||c|c|c|c|c|c|}
  \hline
  	&\multicolumn{2}{c|}{$c_1$ }&\multicolumn{2}{c|}{$c_2$ }&\multicolumn{2}{c|}{$\kappa$ }\\ \cline{2-7}
	Grid&$L_1$-error & EOC & $L_1$-error & EOC & $L_1$-error & EOC\\
	\hline
      {25$\times$25/50$\times$50} & 2.399e-02 & &4.195e-02 & &4.346e-02 & \\
      {50$\times$50/100$\times$100} & 6.067e-03 & 1.9831& 1.088e-02 & 1.9475 & 1.069e-02 & 2.0233\\
      {100$\times$100/200$\times$200} & 1.514e-03 & 2.0026& 2.751e-03 & 1.9831 & 2.681e-03 & 1.9956\\
      {200$\times$200/400$\times$400} & 3.785e-04 & 2.0003& 6.912e-04 & 1.9930 & 6.724e-04 & 1.9954\\
	\hline
    \end{tabular}\end{center}
  \caption{$L_1$-errors and EOC of the components $c_1$, $c_2$, $\kappa$. The parameter set and the initial conditions are described in Experiment \ref{F:exp:conv} (Section \ref{F:sec:exp.desc}). See also Fig. \ref{F:fig:cnv.cpu} (left).} \label{F:tbl:L1EOC}

  \begin{center}\begin{tabular}{|c||c|c|c|c|c|c|}
	\hline
  	&\multicolumn{2}{c|}{$c_1$ }&\multicolumn{2}{c|}{$c_2$ }&\multicolumn{2}{c|}{$\kappa$ }\\ 		\cline{2-7}
	Grid & $L_2$-error & EOC & $L_2$-error & EOC & $L_2$-error & EOC\\
	\hline
      25$\times$25/50$\times$50 & 5.375e-03 & & 1.303e-02 & & 1.309e-02 & \\
      50$\times$50/100$\times$100 & 1.342e-03 & 2.0024 & 3.322e-03 & 1.9713 & 3.187e-03 & 2.0375\\
      100$\times$100/200$\times$200 & 3.350e-04 & 2.0018 & 8.347e-04 &1.9928 & 8.029e-04 & 1.9890\\
      200$\times$200/400$\times$400 & 8.369e-05 & 2.0010 & 2.090e-04 & 1.9977& 2.012e-04 & 1.9965\\
	\hline
    \end{tabular}\end{center}
  \caption{$L_2$-errors and EOC of the components $c_1$, $c_2$, $\kappa$. The setting is described in Experiment~\ref{F:exp:conv}.} \label{F:tbl:L2EOC}

\end{table}

\setlength{\tabcolsep}{12pt}
\begin{figure}[t] 
	\centering
	\resizebox{\linewidth}{!}{
		\begin{tabular}{c c}
			\includegraphics{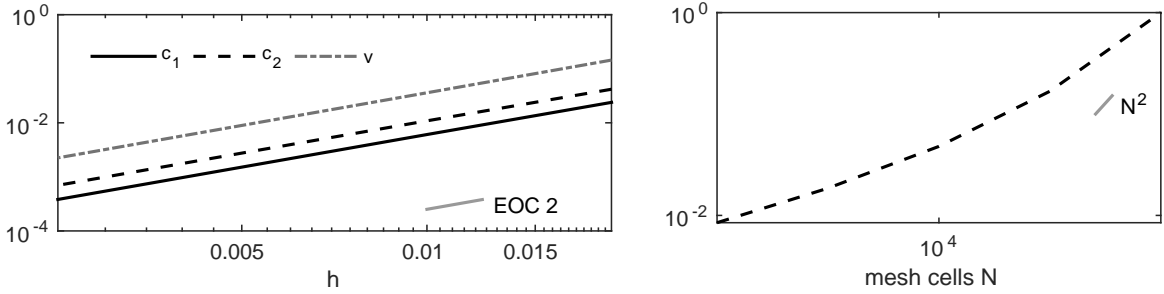} &\includegraphics{cputime} 
		\end{tabular}}
		\vfill \vspace{.25em}
		\caption{(Left:) Graphical representation of the convergence order for the components $c_1$, $c_2$, and $v$ on a two dimensional domain with the grid step size $h$, see also Tables \ref{F:tbl:L1EOC}, \ref{F:tbl:L2EOC}. (Right:) Graph of the relative computational costs as a function of the total number of computational cells. The results correspond to Experiment \ref{F:exp:conv}.}\label{F:fig:cnv.cpu}
\end{figure}

In Tables \ref{F:tbl:L1EOC}, \ref{F:tbl:L2EOC} and in Fig. \ref{F:fig:cnv.cpu} (left) we present the \textit{experimental order of convergence} (EOC) rate for the developed method using the parameters and initial data from the Experiment \ref{F:exp:conv}. We can clearly recognize the second order convergence in all five components of the system. The corresponding computational costs are presented in the Fig. \ref{F:fig:cnv.cpu} (right), where the actual values have been scaled with respect to the more expensive grid $400\times 400$.

\begin{figure}[t]
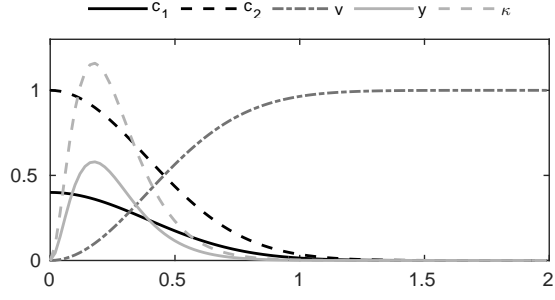
 
	\centering
		\begin{tabular}{b{0.5\linewidth}}
			~\hspace{5.2em}\includegraphics[width=0.695\linewidth]{legendAll}\\
			\hspace{2.6em}\includegraphics[width=0.9\linewidth]{sliceini} \\ 
			\end{tabular}	
			\caption{Radial cut (positive semi-axis) of the two-dimensional initial conditions of the Experiment \ref{F:exp:main}. See also Fig. \ref{F:fig:2D.slice} (right) for the numerical solution at the final time.}
	\label{F:fig:ICs}
\end{figure}

\begin{figure}[t] 
	\centering
	\resizebox{\linewidth}{!}{
		\begin{tabular}{c c}
			& ~\hspace{2.5em} \includegraphics{legendAll}\\
			\includegraphics{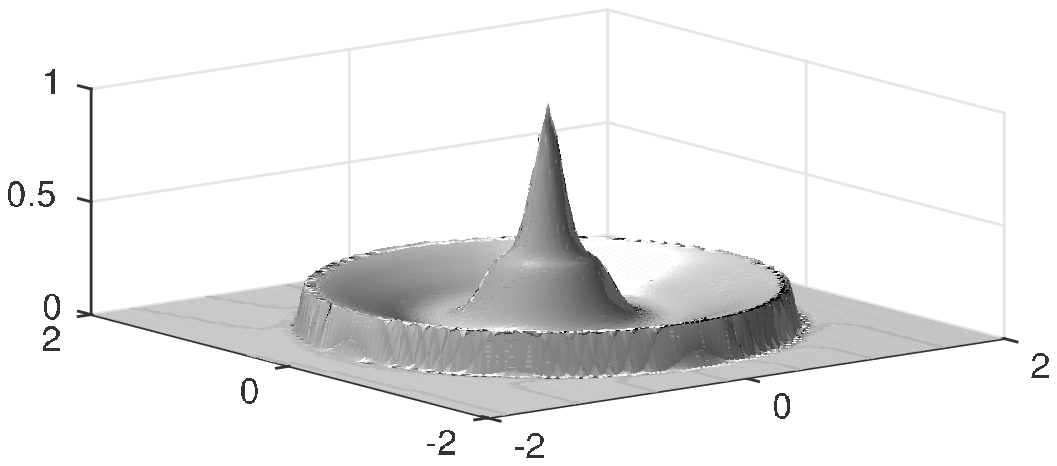} 
				&\includegraphics{slicetau15}
		\end{tabular}}

		\caption{Experiment \ref{F:exp:main}. (Left:) Distribution of the migrating $c_2$ cells for a delay with $\tau=15$ at the time $t=0.5$. The invasion pattern exhibits a steep front not existing in the initial condition, cf. Fig. \ref{F:fig:ICs}. (Right:) Radial cuts of all the components of the solution. The component $c_2$ develops a steep front. The position of the propagating front of $c_2$ and its magnitude  depends on the delay $\tau$, cf. Fig. \ref{F:fig:results}.} \label{F:fig:2D.slice}
\end{figure}

In Figs. \ref{F:fig:ICs} and \ref{F:fig:2D.slice} we see the initial conditions and the final time solutions of (\ref{F:eq:orig2}) according to the the particular Experiment \ref{F:exp:main}. Despite the smoothness of the initial conditions, a steep front is formed in $c_2$, see Fig. \ref{F:fig:2D.slice} (left). All the components of the solution are presented in Fig. \ref{F:fig:2D.slice} (right), where we can also recognize the particular structure of $c_2$ in detail: a propagating front which seems to be ``separated'' from the original part of the tumour is followed by a smooth part.
	
\begin{figure}[t] 
	\centering
	\resizebox{\linewidth}{!}{
		\begin{tabular}{c c}
			\includegraphics{delayAggrStudy}&\includegraphics{delayMassStudy}  \\
			 \vspace{-3em}	$\tau$ & $~~~~~~~~~~~~~~~~~~~~~~\tau$
		\end{tabular}}
		\vfill \vspace{.25em}
		\caption{Presented here is the effect of the delay parameter $\tau$ in the aggressiveness of the tumour. (Left:) The position of the propagation front of the migrating cells decreases with the increase of the delay (solid line). The height of the propagating front of the $c_2$ on the other hand increases with the delay up to $\tau\approx22$ (dashed line). For larger delay value the front height decreases, cf. Fig. \ref{F:fig:slicelargetau}. (Right:) The mass of the migrating cells increases with the delay whereas the mass of the proliferative decreases slightly. The corresponding computational setting is described in Experiment \ref{F:exp:main} (Section \ref{F:sec:exp.desc}).} \label{F:fig:results}
\end{figure}

\begin{figure}[t]
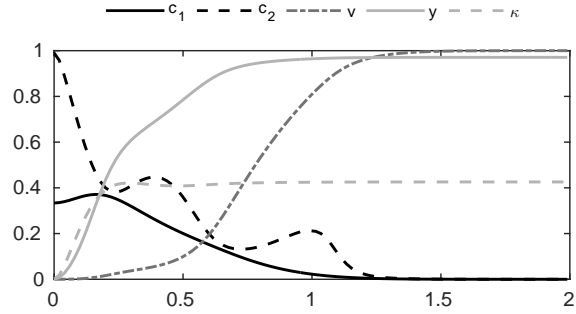
 
	\centering
		\begin{tabular}{b{0.5\linewidth}}
			~\hspace{5.2em}\includegraphics[width=0.68\linewidth]{legendAll}\\
			~\hspace{1.8em}\includegraphics[width=0.93\linewidth]{slicetau30}
		\end{tabular}
	\caption{A radial cut of the numerical solution at the final time of the Experiment \ref{F:exp:main} for a ``large'' delay $\tau=31$. The propagating front has invaded to a lesser extend than for smaller values of $\tau$, cf. Figs. \ref{F:fig:2D.slice} (right) and \ref{F:fig:slicesChiVsTau} (left).}
	\label{F:fig:slicelargetau}
\end{figure}

\begin{figure}[t]
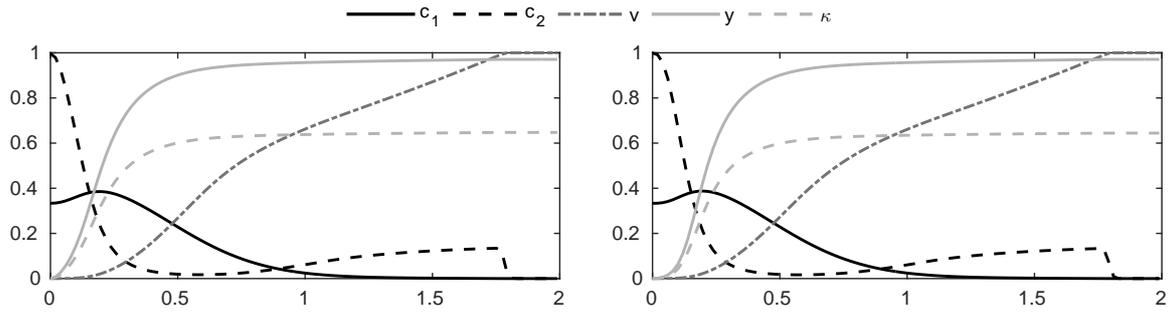
 
	\centering
	\resizebox{\linewidth}{!}{
		\begin{tabular}{c c}
			\multicolumn{2}{c}{\includegraphics[width=25em]{legendAll}}\\
			\includegraphics{slicetauzero} 
			&\includegraphics{sliceMicroChi}  
		\end{tabular}}
		\caption{Graphical comparison between $\tau = 0$ and $\chi=0.001$ (left) versus $\tau=15$ and $\chi =0.0001$ (right). The results at the same final time are almost identical. We can so deduce that the ``convergence'' of the ODE subsystem to the quasi-steady state is very fast, and that the results on the aggressiveness of $c_2$ is mostly due to the compound delay $\chi\tau$ and less due to the time scale $\chi$. The computational setting is described in Experiment \ref{F:exp:main}.}
		\label{F:fig:slicesChiVsTau}
\end{figure}

In Fig. \ref{F:fig:results} we present the dependence of the numerical solution of (\ref{F:eq:orig2}) on the delay parameter $\tau$. In particular, an increasing value of $\tau$ leads to a decreasing position of the propagation front of $c_2$ and an increase of the front height. For values of $\tau\geq 22$ we also see a drop of the front height, which is due to the emerging of a secondary invasion front, see also Fig. \ref{F:fig:slicelargetau}. We can also see a constant rate increase of the mass of $c_2$ with $\tau$ (right), whereas the mass of $c_1$ is not significantly influenced. 

Moreover, we can also see in Fig. \ref{F:fig:slicelargetau} that larger delay values $\tau$ cause the propagating front of the migratory $c_2$ cells to invade to a lesser extend than for smaller $\tau$, cf. Figs. \ref{F:fig:2D.slice} and \ref{F:fig:slicesChiVsTau}. A secondary front, that follows closely the primary front, affects its magnitude, see also Fig. \ref{F:fig:2D.slice}. The experimental setting is given in Experiment \ref{F:exp:main}.

In Fig. \ref{F:fig:slicesChiVsTau} we present a comparison between $\tau = 0$ and $\chi=0.001$ (left) versus $\tau=15$ and $\chi =0.0001$ (right). We note that although the ``stiff'' coefficients in $y$ and $\kappa$ in (\ref{F:eq:orig2}) differ by one order of magnitude (due to $\chi$), the final time results are almost identical. We verify that the ``aggressiveness'' of $c_2$, as we witness in Fig. \ref{F:fig:results}, is mostly influenced by the composite delay $\chi\tau$ and less by the actual time scaling $\chi$.

\subsection{Description of experiments}\label{F:sec:exp.desc}\sectionmark{Description of experiments}

Here we give technical details on the experiments that have been presented in this work. Our simulations have been performed on the computational domain $\Omega =[-2,2]\times[-2,2]$. In all experiments we have employed zero Neumann boundary conditions for the advective-diffusive component $c_2$ of the solution 
\[
	- D_c \frac{\kappa}{1+(c_1+c_2)v}\frac{\partial c_2}{\partial \vec n}+ 
	D_h\frac{ \kappa v}{1+v}c_2\frac{\partial v}{\partial \vec n}=0
\]
where $\vec n$ is the outward normal vector to the computational domain $\Omega$.

\begin{experiment}{0}\label{F:exp:conv}
	This experiment corresponds to the convergence results in Fig. \ref{F:fig:cnv.cpu} and Tables \ref{F:tbl:L1EOC}--\ref{F:tbl:L2EOC}. Following the original derivation of the model \cite{Stinner.2015}, we consider the following set of parameters: 
	\renewcommand{\arraystretch}{1.1}
	\begin{equation}\label{F:exp:conv:params}
		\left\{
		\begin{array}{*4{>{\displaystyle}l}}
		\mu_c=1,		&\qquad \eta_1=0.05, 	&\qquad \gamma=0.055,		&\qquad \lambda=0.076, \\
		D_c=10^{-3},	&\qquad D_h=1, \\
		\delta_v=10, 	&\qquad \mu_v=0.3, 		&\qquad \eta_2=0.9, \\
		k_1=2, 			&\qquad k_{-1}=0.06,\\
		q=3, 			&\qquad M=1,\\
		\chi= 0.01, 	&\qquad \tau=0.04.
		\end{array}
		\right.
	\end{equation}
	The initial condition reads
	\renewcommand{\arraystretch}{1.3}
	\begin{equation}
		\left\{
		\begin{array}{*4{>{\displaystyle}l}}
		c_1(0,x)&=0.4\, e^{-\frac{1}{\varepsilon}\(x^2+y^2\)},\\
		c_2(0,x)&=e^{-\frac{1}{\varepsilon}\(x^2+y^2\)},\\
		v(0,x)&=1-c_2(0,x),\\
		y(0,x)&=20\, f_{\gamma}\(5\(x^2+y^2\),2,15\),\\
		\kappa(0,x)&=2\,y(0,x),
		\end{array}
		\right.
	\end{equation}
	for $\varepsilon=1.5$ and $x \in \Omega$, where we employ the density function
	of the gamma distribution,
	\begin{equation}\label{F:eq:Gamma}
		f_\gamma\(x,a,b \)=\frac{1}{b^a\Gamma(a)}x^{a-1}e^{\frac{-x}{b}},\ \mbox{ where }\ \Gamma(a)=\int_0^\infty t^{a-1}e^{-t} dt.
	\end{equation}
\end{experiment}

\begin{experiment}{1}\label{F:exp:main}
	The parameters and initial conditions that follow, correspond to Figs. \ref{F:fig:ICs}, \ref{F:fig:2D.slice}, \ref{F:fig:results}, \ref{F:fig:slicelargetau}, \ref{F:fig:slicesChiVsTau}. Parameters
	\renewcommand{\arraystretch}{1.1}
	\begin{equation}\label{F:exp:main:params}
		\left\{
		\begin{array}{*4{>{\displaystyle}l}}
		\mu_c=1,		&\qquad \eta_1=0.05, 	&\qquad \gamma=0.055,		&\qquad \lambda=0.152, \\
		D_c=10^{-2},	&\qquad D_h=10, \\
		\delta_v=5, 	&\qquad \mu_v=0.3, 		&\qquad \eta_2=0.9, \\
		k_1=2, 			&\qquad k_{-1}=0.06,\\
		q=3, 			&\qquad M=1,\\
		\chi= 0.01, 	&\qquad \tau=0.04.
		\end{array}
		\right.
	\end{equation}
	Initial conditions
	\renewcommand{\arraystretch}{1.3}
	\begin{equation}
		\left\{
		\begin{array}{*4{>{\displaystyle}l}}
		c_1(0,x)&=0.4\, e^{-\frac{1}{\varepsilon}\(x^2+y^2\)},\\
		c_2(0,x)&=e^{-\frac{1}{\varepsilon}\(x^2+y^2\)},\\
		v(0,x)&=1-c_2(0,x),\\
		y(0,x)&=15\, f_{\gamma}\(80\sqrt{x^2+y^2},3,7\),\\
		\kappa(0,x)&=2\,y(0,x),
		\end{array}
		\right.
	\end{equation}
	where $\varepsilon=1.5$, $x \in \Omega$, and $f_\gamma$ is defined in (\ref{F:eq:Gamma}).
\end{experiment}

\section{Conclusions} \sectionmark{Conclusions}

Since their first derivation cancer growth models have been a theater for the development of new problem-suited numerical methods. This is not only due to the importance of the topic, but more importantly it is also due to complex dynamics of the solutions. Our work aims to be a contribution along these lines. 

We solve numerically the model (\ref{F:eq:orig2}) that was proposed in \cite{Stinner.2015}. The method we employ is a concatenation of a robust, positivity preserving FV method in space with a third order IMEX method in time. The additional challenges that we have addressed are the non-constant diffusion coefficients in $c_2$, the time delay in $\kappa$, as well as the stiff reaction terms in $y$ and $\kappa$. 

We have discretized the non-constant diffusion coefficient using central differences and solved the (implicit) linear system by the Krylov method. For the delay term we perform an interpolation in time between a small number of previously ``saved'' time steps. We treat the stiffness applying a secondary condition (besides the CFL) by adjusting the time step of the method. The additional condition leads to an adaptive time stepping by employing an explicit Euler step of the $\kappa$ equation.

We verify numerically that our method is second order accurate and identify its computational cost. Our extensive numerical experiments indicate that the migrating cancer cells develop a steep propagating front. We also show that its aggressiveness depends on the time delay.

\paragraph{Acknowledgements:}
	The authors wish to thank Christina Surulescu and Christian Stinner for the fruitfull discussions during the preparation of this work.

%
\bibliographystyle{plain}
	\bibliography{Kol-Luk-Sfak-Wie}

\end{document}